\documentclass[12pt]{article}

\usepackage{amssymb,latexsym   }
\usepackage[centertags]{amsmath}

\numberwithin{equation}{section}

\renewcommand{\geq}{\geqslant}
\renewcommand{\leq}{\leqslant}

\renewcommand{\ge}{\geqslant}
\renewcommand{\le}{\leqslant}

\newtheorem{thm}[equation]{Theorem}
\newtheorem{lemma}[equation]{Lemma}

\newtheorem{definition}[equation]{Definition}

\newcommand{\N}{\mathbb{N}}

\newcommand{\R}{\mathbb{R}}
\newcommand{\C}{\mathbb{C}}
\newcommand{\D}{\mathbb{D}}

\newcommand{\A}{\operatorname{Area}}
\newcommand{\T}{\mathbb{T}}
\newcommand{\la}{\lambda}

\newcommand{\ga}{\gamma}

\DeclareMathOperator*{\osc}{osc}

\newcommand{\Z}{\mathbb{Z}}

\renewcommand{\arg}{{\operatorname{arg}}}
\newcommand{\cN}{{\mathcal{N}}}
\newcommand{\cF}{{\mathcal{F}}}
\newcommand{\Area}{{\operatorname{Area}}}

\begin{document}

\title{Sign and area in nodal geometry \\ of Laplace eigenfunctions}

\author{F\"edor Nazarov, Leonid Polterovich and Mikhail Sodin}

\maketitle

\begin{abstract}
The paper deals with asymptotic nodal geometry for the
Laplace-Beltrami operator on closed surfaces. Given an
eigenfunction  $f$ corresponding to a large eigenvalue, we study
local asymmetry of the distribution of $\text{sign}(f)$ with
respect to the surface area. It is measured as follows:  take any
disc centered at the nodal line $\{f=0\}$, and pick at random a
point in this disc. What is the probability that the function
assumes a positive value at the chosen point? We show that this
quantity may decay logarithmically as the eigenvalue goes to
infinity, but never faster than that. In other words, only a mild
local asymmetry may appear. The proof combines methods due to
Donnelly-Fefferman and Nadirashvili with a new result on harmonic
functions in the unit disc.
\end{abstract}

\section{Introduction and main results}\label{sect1}

Consider a compact manifold $S$ endowed with a $C^{\infty}$
Riemannian metric $g$. Let $\{f_{\lambda}\}$, $\lambda \nearrow
+\infty$, be any sequence of eigenfunctions of the
Laplace-Beltrami operator $\Delta_g$:
$$
\Delta_g f_{\lambda} + \lambda f_{\lambda}=0\,.
$$
The eigenfunctions $f_{\lambda}$ give rise to an interesting
geometric object, nodal sets $L_{\lambda} = \{f_{\lambda}=0\}$.
Each $L_{\lambda}$ is a closed hypersurface with quite tame
singularities. For instance, when $S$ is 2-dimensional, any nodal
line $L_{\lambda}$ at a singular point $p$ looks like the union of
an even number of smooth rays  meeting at $p$ at equal angles
\cite[Chapter~III]{SY}. In spite of this ``infinitesimal
simplicity", the global picture of nodal sets for large $\lambda$
becomes more and more complicated. This is partially due to the
fact that $L_{\lambda}$ is $\sim 1/\sqrt{\lambda}$-dense in $S$.

Asymptotic geometry of nodal sets as $\lambda \nearrow \infty$
attracted a lot of attention of both mathematicians and physicists
though it is still far from being understood  (see \cite{JNT, BGS}
for discussion on recent developments). The idea of studying the
asymptotic behavior comes from quantum mechanics where
$f_{\lambda}^2$ (properly normalized) is interpreted as the
probability density of the coordinate of a free particle in the
pure state corresponding to $f_{\lambda}$, and $\lambda \nearrow
+\infty$ corresponds to the quasi-classical limit.

A nodal domain is a connected component of the set $S \setminus
L_{\lambda}$. All nodal domains can be naturally grouped into two
subsets $S_+(\lambda):= \{f_{\lambda}> 0\}$ and $S_-(\lambda):=
\{f_{\lambda}< 0\}$. Our story starts with two fundamental results
obtained by Donnelly and Fefferman.

The first one is {\it a ``local version" of the Courant nodal
domain theorem} \cite{DF3}: let $D \subset S$ be a metric ball and
let $U$ be any component of $S_+(\la) \cap D$  such that
\begin{equation} \label{eqdeep}
U \cap \frac{1}{2}D \neq \varnothing\,,
\end{equation}
which means that $U$ enters deeply enough into $D$ \footnote {Here
and below $\frac{1}{2}D$ stands for the ball with the same center
as $D$ whose radius equals half of the radius of $D$. The radii of
all metric balls are assumed to be less than the injectivity
radius of the metric. }. Then
\begin{equation}
\label{eqDF} \frac{\text{Volume}(U)}{\text{Volume}(D)} \geq a\cdot
\lambda^{-k}
\end{equation}
where $a$ depends only on the metric $g$ and $k$ only on the
dimension of $S$. The rate of decay of the right hand side (a
negative power of $\lambda$) cannot be improved --- a suitable
example can be easily produced already in the case of standard
spherical harmonics. The sharp value of the constant $k$ is,
however, still unknown (see papers \cite{ChM, Lu} for estimates on
$k$).

The second result is the following {\it quasi-symmetry theorem}
proved in \cite[p.~182]{DF1} under the extra assumption that the
metric $g$ is real analytic. Let $D \subset S$ be a fixed ball.
Then there exists $\Lambda $ depending on the radius of the ball
$D$ and the metric $g$ such that for all $\lambda > \Lambda$
\begin{equation}
\label{eqDF1} \frac{\text{Volume}(S_+(\lambda) \cap
D)}{\text{Volume}(D)} \geq a,
\end{equation}
where $a>0$ depends only on the metric $g$.

From the geometric viewpoint, there is a significant difference
between the measurements presented above: the quasi-symmetry
theorem \eqref{eqDF1}  deals with a ball of fixed radius and large
$\lambda$. In contrast to this, the local version of the Courant
theorem \eqref{eqDF} is valid for all scales and all $\lambda$'s
though the collection of balls depends on $\lambda$ through the
``deepness assumption" \eqref{eqdeep}. A natural problem arising
from this discussion is to explore what remains of quasi-symmetry
on all scales and for all $\lambda$, provided that the nodal set
enters deeply enough into a ball: $L_{\lambda} \cap \frac{1}{2}D
\neq \varnothing$.

In the present paper we deal with this problem in the case when
$S$ is a compact connected surface and the metric $g$ is
$C^{\infty}$-smooth. Our main finding is that only a mild local
asymmetry may appear. If in formula \eqref{eqDF} one replaces a
single component $U$ by the whole set $S_+(\lambda)$, the right
hand side changes its behavior: instead of a negative power of
$\lambda$, it becomes $(\log \lambda)^{-1}
(\log\log\lambda)^{-1/2}$. Moreover, we will show that even for
the standard spherical metric it cannot be better than
$(\log{\lambda})^{-1}$. We believe that the double logarithm
factor reflects a deficiency in our method (see discussion in
Section~\ref{subsec7.1}). The precise formulations follow.

\begin{thm}\label{thmsymmet}
Let $S$ be a compact connected surface endowed with a smooth
Riemannian metric $g$, and let $f_{\lambda},\; \lambda \geq 3$, be
an eigenfunction of the Laplace-Beltrami operator. Assume that the
set $S_+(\lambda):= \{f_{\lambda}> 0\}$ intersects a metric disc
$\frac{1}{2}D$. Then
$$
\frac{\A(S_+(\lambda) \cap D)}{\A(D)} \geq \frac{a}{\log \lambda
\cdot \sqrt{\log\log\lambda} }
$$
where the constant $a>0$ depends only on $g$.
\end{thm}
(The condition $\lambda\geq 3$ is imposed here only because
$\log{}$ and $\sqrt{}$ are not defined if the argument is less
than $0$. For $\lambda<3$ the theorem holds with the right hand
side replaced by a constant $a>0$ depending on $g$ only.)

The next result illustrates sharpness of the previous estimate up
to the double logarithm \footnote{It is worth mentioning that on
``microscopic scales", when the radii of discs $D$ are less than $
(\la \log\log \la)^{-1/2} $, our approach gives an optimal bound
$$
\frac{\A(S_+(\la) \cap D)}{\A (D)} \geq \frac{a}{\log\la}\,.
$$}:
\begin{thm}\label{thmsharp}
Consider the 2-sphere $\mathbb{S}^2$ endowed with the standard
metric. There exist a positive numerical constant $C$, a sequence
of Laplace-Beltrami eigenfunctions $f_i,\; i \in \N$ corresponding
to eigenvalues $\lambda_i \to \infty$, and a sequence of discs
$D_i \subset \mathbb{S}^2$ such that each $f_i$ vanishes at the
center of $D_i$ and
$$
\frac{\A(S_+(\lambda_i) \cap D_i)}{\A(D_i)} \leq \frac{C}{\log
\lambda_i}\,.
$$
\end{thm}

After Donnelli and Fefferman \cite{DF1}, various versions of
quasi-symmetry for eigenfunctions were studied by Nadirashvili
\cite{N} and Jakobson-Nadirashvili \cite{JN}. To a high extent,
the present research was stimulated by Nadirashvili's article
\cite{N}.

Our approach to Theorem~\ref{thmsymmet} is based on the analysis
of the eigenfunctions $f_{\lambda}$ on discs of radius $\sim
1/\sqrt{\lambda}$. The proof consists of four main ingredients
that we are going to describe right now. The first three of them
exist in the literature. Our innovation is the last one, namely,
the calculation of the asymptotical behaviour of the Nadirashvili
constant for harmonic functions.

\medskip\noindent{\sc Donnelly-Fefferman growth bound}.
For any continuous function $f$ on a closed disc $D$ (in any
metric space), define its {\em doubling exponent} $\beta(D, f)$ by
$$
\beta(D, f) = \log\frac{\max_{D} |f|}{\max_{\frac12D} |f|}.
$$
The following fundamental inequality was established in \cite{DF1}
in any dimension. For any metric disc $D \subset S$ and any
$\lambda$,
\begin{equation}\label{eqdoubl}
\beta(D, f_\lambda) \leq a \sqrt{\lambda}
\end{equation}
where the constant $a$ depends only on the metric $g$.

\medskip\noindent{\sc Reduction to harmonic functions.}
Assume now that $D \subset S$ is a disc of radius  $\sim
1/\sqrt{\lambda}$. It turns out that on this scale the
eigenfunction $f_{\lambda}$ can be ``approximated" by a harmonic
function $u$ on the unit disc $\D$. More precisely, the set
$\{f_{\lambda}>0\}$ can be transformed into the set $\{u>0\}$ by a
$K$-quasiconformal homeomorphism with a controlled dilation $K$.
Moreover, the doubling exponent of $u$ on $\D$ is essentially the
same as that of $f_\lambda$ in $D$. This idea is borrowed from
Nadirashvili's paper \cite{N}. The details are presented in
Section~\ref{sect3} below.

\medskip\noindent{\sc Topological interpretation of the doubling
exponent.} Let \\
$u\colon\D\to\R$ be a non-zero harmonic function. Denote by
$\nu(r\T, u)$ the number of sign changes of $u$ on the circle $r\T
= \{|z|=r\}$. Then
\begin{equation}\label{eqgelfond}
C^{-1} ( \beta( \tfrac14 \D, u) - 1) \leq \nu( \tfrac12 \T, u)
\leq C (\beta(\D, u) + 1)
\end{equation}
where $C$ is a positive numerical constant. This result goes back
to Gelfond~\cite{Gelfond} (cf. \cite{Robertson}, \cite{HK}, and
\cite[Theorem~3]{KY}). We will need the inequality on the right
only, which will be proved in Section~\ref{sect2}. The inequality
on the left is presented here just for completeness.

\medskip\noindent{\sc The Nadirashvili constant.}
Denote by ${\cal H}_d$ the class of all non-zero harmonic
functions $u$ on $\D$ with $u(0) = 0$ that have no more than $d$
sign changes on the unit circle $\T$. Define the {\em
Nadirashvili constant}
$$
\mathcal N_d := \inf_{u \in {\cal H}_d} \A(\{u>0\})\,.
$$
Using an ingenious compactness argument, Nadirashvili \cite{N}
showed that $\mathcal N_d$ is strictly positive. Our next result
gives a satisfactory estimate of the Nadirashvili constant:
\begin{thm}\label{thmnadconst} There exists a positive numerical constant $C$
such that for each $d\ge 2$,
$$
\frac{C^{-1}}{\log d} \leq \mathcal N_d \leq \frac{C}{\log d}\,.
$$
\end{thm}

Nadirashvili's proof of positivity of $\mathcal N_d$ is
non-constructive, hence we had to take a different route. Our
approach is based on one-dimensional complex analysis.

\medskip
The four steps described above yield Theorem~\ref{thmsymmet} in
the case when the disc $D$ is {\em small}, that is, of radius
$\leq a \lambda^{-1/2}$. The double logarithm term is the price we
pay for the fact that the transition from the eigenfunction
$f_{\lambda}$ to the approximating harmonic function $u$ is given
by a quasiconformal homeomorphism, which in general is only
H\"older. The case of an arbitrary (not necessarily small) disc
$D$ is based on the following standard argument. The nodal line
$L=\{f_{\lambda} = 0\}$ is $\sim 1/\sqrt{\lambda}$-dense in $S$
(see e.g. \cite{Bruning}). Hence every disc $D$ with $L \cap
\frac{1}{2} D \neq \varnothing$ contains a disjoint union of small
discs $D_i$ whose centers lie on $L$ and such that the total area
of these discs is $\geq \text{const}\cdot\A(D)$. Since the area
bound is already established for each $D_i$, it extends with a
weaker constant to $D$. This completes the outline of the proof of
Theorem \ref{thmsymmet}.

\medskip\par\noindent {\sc Organization of the paper. }
The next section is devoted to harmonic functions on the unit
disc. We establish the lower bound $\mathcal N_d \geq c(\log
d)^{-1}$ for the Nadirashvili constant and prove the right
inequality in \eqref{eqgelfond} relating the number of boundary
sign changes to the doubling exponent.

In Section~\ref{sect3}, we deal with solutions of the
Schr\"odinger equation in the unit disc with small potential. This
Schr\"odinger equation is nothing else but an appropriately
rescaled equation $\Delta_g f + \lambda f = 0$ written in local
conformal coordinates on the surface. For the solutions $F$ of
this equation, we prove a lower bound on $\A(\{F>0\}$ in terms of
the doubling exponent of $F$. The proof is based on a
quasiconformal change of variables that reduces the problem to the
estimate for harmonic functions obtained in Section \ref{sect2}.

In Section \ref{sect4}, we present the easiest proof of the
Donnelli-Fefferman fundamental inequality \eqref{eqdoubl} we are
aware of. We should warn the reader that our proof works in
dimension $2$ only. It is based on a simple observation about
second order linear ODEs in Hilbert spaces. The reader familiar
with the Donnelli-Fefferman inequality \cite{DF1, JL} can
disregard this section.

At this point we have all the ingredients necessary to prove
Theorem~\ref{thmsymmet}. This is done in Section~\ref{sect5}.

In Section \ref{sect6}, we present examples illustrating the local
logarithmic asymmetry for harmonic functions and Laplace-Beltrami
eigenfunctions. Our construction uses the complex double
exponential function $\exp \exp z$. We confirm the upper bound for
the Nadirashvili constant $\mathcal N_d$, which completes the
proof of Theorem~\ref{thmnadconst}. Then, ``transplanting" the
Taylor series of the obtained harmonic function at $0$ to the
north pole of the unit sphere, we obtain Theorem~\ref{thmsharp}
that shows that our main result is already sharp for spherical
harmonics up to the double logarithm.

The paper concludes with discussion and questions. In particular,
we indicate a link between the expectation of the doubling
exponent of an eigenfunction $f_{\lambda}$ on a random disc of
radius $\sim 1/\sqrt{\lambda}$ and the length of its nodal line
$\{f_{\lambda} = 0\}$.

\medskip\par\noindent{\sc Convention. }
Throughout the paper, we denote by $c$, $c_0$, $c_1$, $c_2$, ...
positive numerical constants, and by $a$, $a_0$, $a_1$, ...
positive constants that depend only on the metric $g$. In each
section we start a new enumeration of these constants.

\subsection*{Acknowledgement} We thank Kari Astala,
Eero Saksman, Uzy Smilansky, and Sasha Volberg for helpful
discussions and comments.


\section{The area estimate for harmonic functions}\label{sect2}

In this section, we show that for any non-zero harmonic function
$u$ on $\D$ vanishing at the origin,
\begin{equation}\label{eq5.1}
\text{Area}(\{u > 0\}) \geq \frac{c}{\log\nu (\T, u)}\,,
\end{equation}
i.e., we prove the lower bound for the Nadirashvili constant
$\mathcal N_d$ in Theorem~\ref{thmnadconst}. Then we prove the
right hand part of estimate \eqref{eqgelfond}. Together with
\eqref{eq5.1}, it yields
\begin{thm}\label{thmharm} Let $u$ be a non-zero harmonic function on the
unit disc $\D$ vanishing at the origin. Then
$$
\A(\{u>0\}) \ge \frac{c_0}{\log\beta^*(\D, u)}
$$
where $\beta^*:=\max(\beta, 3)$.
\end{thm}

Consider the analytic function $f\colon \D \to \C$ with
$\operatorname{Re}f = u$ and $f(0) = 0$. Assume that $f$ does not
vanish on $r\T$. Consider all arcs $L\subset r\T$ travelled
counterclockwise (including the entire circumference $r\T$ viewed
as an arc whose end and beginning coincide). Put
$$
\omega(r\T, f) :=\max_{L\subset r\T}\Delta_L\arg f
$$
where $\Delta_L\arg f$ is the increment of the argument of $f$
over $L$, that is,
$$
\Delta_L \arg f = \arg f(\theta_2) - \arg f(\theta_1)\,,
$$
for $L=[\theta_1; \theta_2]$.
 We shall prove
\begin{thm}\label{thmosc}
Let $f$ be an analytic function on $\D$ vanishing at the origin.
Assume that $f|_{\T}\ne 0$. Then
\begin{equation}\label{eqosc}
\Area(\{\operatorname{Re} f>0\}) \ge \frac{c_1}{\log \omega(\T,
f)}\,.
\end{equation}
\end{thm}

Since $\omega (\T, f) \leq \pi (\nu (\T, \hbox{Re} f)+1)$, this
yields estimate \eqref{eq5.1} and, therefore, the lower bound for
the Nadirashvili's constant.

\medskip\par\noindent{\em Proof of Theorem~\ref{thmosc}: }
For $k\in\mathbb N$, denote by $\cF_k$ the class of analytic
functions $f$ on $\D$ such that $f(0) = 0$, $f$ does not vanish on
$\T$, and $\omega(\T, f) \leq 2\pi\cdot 2^k$. Put
$$
A_k = \inf_{f \in \cF _k} \Area(\{\operatorname{Re} f > 0\})\,.
$$
The estimate \eqref{eqosc} would follow from the inequality
\begin{equation}\label{eqind}
A_k \geq \frac{c_2}{k}.
\end{equation}

Start with any $f \in \cF_ k$, and define $\delta$ by
$$
1-2\delta = \sup\{r\colon f|_{r\T} \ne 0, \ \omega(r\T, f) <
2\pi\cdot 2^{k-1}\}\,.
$$
If this set is empty, we simply take $\delta = \frac12$.

Consider the annulus $E=\{1-2\delta<|z|<1-\delta\}$ and its subset
$E_+ = \{z\in E\colon \operatorname{Re}f(z)> 0 \}$. The heart of
our argument is the following
\begin{lemma}\label{lemmafund}
$\A(E_+) \geq c_3 \delta^2$.
\end{lemma}

Assuming the lemma, let us prove inequality \eqref{eqind} by
induction on $k$. First of all, consider the case $k=1$. Since $f$
vanishes at the origin,
$$
\omega (r\T, f) \ge \Delta_{r\T} \arg f \ge 2\pi
$$
for all $r>0$. Therefore, we can take $\delta = \frac12$, and
Lemma \ref{lemmafund} yields $A_1 \geq \frac14 c_3$. Hence, taking
$c_2= \frac14 c_3$, we prove the induction base for claim
\eqref{eqind}.

Assume now that \eqref{eqind} is true for $k-1$. Let us prove it
for $k$. Take any $f\in \cF_k$. If $\delta=\frac{1}{2}$, Lemma
\ref{lemmafund} immediately yields
$$
\A(\{\operatorname{Re}f > 0\}) \geq \text{Area}(E_+)\geq c_2\geq
\frac{c_2}{k}\,.
$$
Otherwise, we can find $r>0$ arbitrarily close to $1-2\delta$ and
such that $f$ does not vanish on $r\T$ and $\omega(r\T, f) <
2\pi\cdot 2^{k-1}$. Put $g(z) = f(rz)$, $z \in \D$. Note that $g
\in \cF_{k-1}$ due to our choice of $r$. Obviously,
$$
\A(\{\operatorname{Re}f > 0\}) \geq \text{Area}(E_+) + r^2
\text{Area}(\{\operatorname{Re}g > 0\})\,.
$$
Applying Lemma \ref{lemmafund} and the induction assumption and
letting $r\to 1-2\delta$, we get
$$
A_k \geq c_3 \delta^2 + (1-2\delta)^2 \frac{c_2}{k-1} =
\frac{c_3}{4} (2\delta)^2 + (1-2\delta)^2 \frac{c_2}{k-1} \,.
$$
Note that the minimal value of the function $q(x)= \alpha x^2 +
\beta (1-x)^2$ equals $\alpha\cdot \beta /(\alpha + \beta)$. Thus,
$$
A_k \ge \frac{(c_3/4) \cdot (c_2/(k-1))}{c_3/4 + c_2/(k-1) } =
\frac{c_2}{k + 4c_2/c_3 -1}\,.
$$
Hence, making the same choice $c_2= c_3/4$ as above, we get $A_k
\geq c_2/k$, and inequality \eqref{eqind} follows. This yields
Theorem~\ref{thmosc} modulo Lemma~\ref{lemmafund}. \hfill $\Box$

\medskip\par\noindent{\em Proof of Lemma~\ref{lemmafund}: }
The proof is based on comparing the upper and the lower bounds for
the integral
$$
\iint_{E_+} |\nabla\arg f|\, d\A\,.
$$
Any function $f\in \cF_k$ admits the factorization
$$
f(z)=e^{g(z)}\prod_{\zeta\in \cN(f)}(z-\zeta)
$$
where $\cN(f)$ is the set of zeroes of $f$ in $\D$ counted with
their multiplicities and $g$ is an analytic function in $\D$. Put
$M:=2^k \cdot 2\pi$. Applying the argument principle, we conclude
that the number $N$ of zeroes of $f$ in $\D$ satisfies
\begin{equation}
\label{eq4.2} N \leq \frac{M}{2\pi}\,.
\end{equation}
Further, for $|\zeta| <1$, the function $\theta \to
\text{arg}(e^{i\theta}-\zeta)$ increases with $\theta$. Therefore,
considering the arc $L\subset \T$ joining the point of the minimum
of $\operatorname{Im}g$ to the point of the maximum of
$\operatorname{Im}g$ counterclockwise, we obtain
\begin{equation}
\label{eq4.1} \osc_{\mathbb T} \operatorname{Im}g := \max_{\mathbb
T}\operatorname{Im} g -\min_{\mathbb T}\operatorname{Im} g
\leq\Delta_L\arg f \leq M.
\end{equation}

Fix $r \in (1-2\delta,1-\delta)$ such that $r\T \cap \cN(f) =
\varnothing$. We call an open arc $I \subset r\mathbb T$ {\it a
traversing arc} if its image curve $f(I)$ traverses the right
half-plane, that is, a continuous branch of $\arg f$ maps $I$ onto
an interval $J=(-\frac{\pi}2 +2\pi m; \frac{\pi}2+ 2\pi m)$ for
some $m\in \Z$. Each traversing arc lies in the set $E_+$ which we
are studying. By our choice of $\delta$, the increment of the
argument of $f$ over some arc $L\subset r \mathbb T$ is at least
$M/2$. Hence $L\cap E_+$ (and, thereby, $r\T\cap E_+$) contains
either at least $M/(4\pi)$ pairwise disjoint traversing arcs or
$\frac{M}{4\pi}-1$ traversing arcs and two ``tails". These tails,
taken together, are as good for our purposes as one full
traversing arc.

Given a traversing arc $I \subset r\mathbb T$, note that
$$
\int_I |\nabla\text{arg} f(z)|\,|dz| \geq \pi\,.
$$
Summing up these inequalities over all traversing arcs lying on
$r\T$ and integrating over $r \in (1-2\delta; 1-\delta)$, we  get
\begin{equation}
\label{eq4.3} \iint_{E_+}|\nabla\text{arg} f|\, d\A \geq
\frac{M\delta}{4}\,.
\end{equation}

On the other hand,
\begin{equation}\label{eqlog_der}
|\nabla\text{arg} f (z)| \leq |\nabla \operatorname{Im}g (z)| +
\sum_{\zeta\in \cN(f)}\frac{1}{|z-\zeta|}\,.
\end{equation}

Next, we use an estimate for the gradient of a harmonic function
$v$ in a disc $D$ of radius $t$ centered at $c$:
$$
|\nabla v(c)| \le \frac2{t} \max_{\partial D} |v|\,,
$$
which easily follows by differentiation of the Poisson integral
representation for $v$ in $D$. Applying this estimate to
$v=\operatorname{Im}g - m$ with $m = \frac12 (\max_\T \mbox{Im} g
+ \min_\T \mbox{Im} g) $ and taking into account inequality
\eqref{eq4.1}, we readily get that
\begin{equation}\label{eqnabla_g}
|\nabla \operatorname{Im}g| \le \delta^{-1} \osc_{\T}
\operatorname{Im} g  \le \frac{M}{\delta}
\end{equation}
everywhere in the annulus $E$. Further,
\begin{eqnarray}
\label{eq4.5} \iint_{E_+}\frac{d\A(z)}{|z-\zeta|}
&=& \iint_{\zeta + E_+} \frac{d\A (w)}{|w|} \nonumber \\
&\le& \iint\limits_{|w|\le \sqrt{\A(E_+)/\pi} } \frac{d\A(w)}{|w|}
\leq 2\sqrt{\pi \A(E_+)}\,. \qquad
\end{eqnarray}

Estimates \eqref{eqlog_der}, \eqref{eqnabla_g}, \eqref{eq4.5}, and
\eqref{eq4.2} give us
$$
\iint_{E_+} |\nabla \arg f|\, d\A \le \frac{M}{\delta} \cdot
\Area(E_+) + \frac{M}{2\pi} \cdot 2\sqrt{\pi \text{Area}(E_+)}\,.
$$
Juxtaposing this  with \eqref{eq4.3} and canceling the factor $M$,
we get
$$
\frac{\delta}{4} \le \frac{\Area(E_+)}{\delta} +
\sqrt{\Area(E_+)/\pi}\,.
$$
This yields $\text{Area}(E_+) \ge c_3 \cdot \delta^2$, proving the
lemma. \hfill $\Box$

\medskip
In order to get Theorem~\ref{thmharm}, we need the following

\begin{lemma}\label{lemmahk} Let $u$ be a non-zero harmonic
function on $\D$ vanishing at the origin. Then $\nu( \frac12 \T,
u) \le c_4 \beta^* (\D, u)$.
\end{lemma}

\medskip\par\noindent{\em Proof of Lemma~\ref{lemmahk}: } The proof is a
minor variation of the argument used in \cite{Gelfond, HK}.
Consider the function
\begin{equation}\label{U}
U(\theta) = u(\tfrac12 e^{i\theta}) = \sum_{k\in\mathbb Z}
\widehat{u}(k) 2^{-|k|} e^{ik\theta}\,,
\end{equation}
where $\{\widehat{u}(k)\}$ are the Fourier coefficients of the
function $\theta\mapsto u(e^{i\theta})$. Since
$|\widehat{u}(k)|\le \max_\D |u|$, we see by inspection of formula
\eqref{U} that the function $ U $ has an analytic extension onto
the strip $ \Pi = \{|\operatorname{Im}\theta| \le \log\sqrt2 \}$
and that
\begin{equation}\label{upper}
\max_{\theta\in\Pi} U(\theta) \le \left(\sum_{k\in\mathbb Z}
2^{-|k|/2}\right) \cdot \max_{\D}|u| = c_5 \cdot \max_{\D} |u|\,.
\end{equation}
At the same time,
\begin{equation}\label{lower}
\max_{\theta\in\mathbb R} |U(\theta)| = \max_{0.5\D} |u| =
e^{-\beta (\D, u)} \cdot \max_{\D} |u|\,.
\end{equation}
Now observe that $\nu (1/2\T, u) $ does not exceed the number of
zeroes of $U$ on the interval $[-\pi, \pi]$. The latter can be
easily estimated using Jensen's formula.

For this purpose, consider the rectangle $ P = \{ |x|\le
\frac{3\pi}{2}, |y|\le \log\sqrt2 \}$ and a conformal mapping $h:
\D \to P$ with $h(0)$ chosen in such a way that
$$
\max_{[-\pi, \pi]} |U| = |(U\circ h)(0)|\,.
$$
There exists $c_6<1$ such that, for all such mappings $h$, we have
$h^{-1}[-\pi, \pi]\subseteq \{|z|\le c_6\}$. Denote by $n(t)$ the
number of zeroes of the analytic function $U\circ h$ in the closed
disc $\{|z|\le t\}$. Then
\begin{eqnarray*}
\nu(1/2\T, u) &\le& n(c_6) \le \frac1{\log (1/c_6)}\int_{c_6}^1
\frac{n(t)}{t}\, dt \le
\frac1{\log(1/c_6)}\int_0^1 \frac{n(t)}{t}\, dt \\ \\
&\stackrel{\mbox{`Jensen'}}=& \frac1{\log(1/c_6)}\, \left(
\frac1{2\pi}\int_{-\pi}^\pi
\log|(U\circ h)(e^{i\theta})|d\theta - \log|(U\circ h)(0)| \right) \\ \\
&\le& \frac1{\log(1/c_6)} \left(\log\max_{\Pi} |U| - \log\max_{\R}
|U|\right)\,.
\end{eqnarray*}
Now, taking into account estimates \eqref{upper} and
\eqref{lower}, we readily see that the right hand side is $ \le
(\log (1/c_6))^{-1}( \log c_5 + \beta(\D, u))$, proving the lemma.
\hfill $\Box$



\section{An area estimate for solutions to Schr\"odinger's
equation with small potential}\label{sect3}

In local conformal coordinates on the surface $S$, the equation
$\Delta_g f_\la + \la f_\la = 0$ reduces to $\Delta f + \la q f =
0$. If the size of the local chart is comparable to the wavelength
$\la^{-1/2}$, then, after rescaling and absorbing the spectral
parameter $\la$ into the potential $q$, one arrives at the
Schr\"odinger equation
\begin{equation}\label{eq1}
\Delta F + q F = 0
\end{equation}
with a bounded smooth potential $q$ on the unit disc $\D$. The
disc is endowed with the complex coordinate $z=x+iy$.

Throughout this section, we assume that $\|q\|:=\max_\D |q| <
\varepsilon_0$ where $\varepsilon_0$ is a sufficiently small
positive numerical constant. The result of the present section is
an intermediate step between Theorems~\ref{thmharm} and
\ref{thmsymmet}.
\begin{thm}\label{thmschro}
Let $F$ be any non-zero solution of equation \eqref{eq1} with
$F(0)=0$. Set
$$
\beta (F) := \sup_{D\subset \D} \beta (D, F)
$$
(the supremum is taken over all discs $D\subset \D$). Set
$\beta^*(F) = \max(\beta (F), 3)$. Then
$$
\A(\{F>0\}) \ge
\frac{c}{\log\beta^*(F)\cdot\sqrt{\log\log\beta^*(F)} }.
$$
\end{thm}

The proof is based on Theorem~\ref{thmharm} and on a chain of
lemmas. By $\|\cdot\|$ we always mean the uniform norm in $\D$.

\begin{lemma}\label{lemma1} If $\varepsilon_0$ is sufficiently small, then
equation \eqref{eq1} admits a positive solution $\varphi$ with
$$
1-c_1\|q\| \le \varphi \le 1\,.
$$
\end{lemma}
\par\noindent{\em Proof of Lemma~\ref{lemma1}: }
Define recursively a sequence of functions $F_i$ by $F_0=1$;
$\Delta F_{i+1} = -qF_i$, $F_{i+1}|_\T = 0$. Then $F_{i+1}$ can be
represented in $\D$ as Green's potential of the function $qF_i$:
$$
F_{i+1}(z) = \iint_\D \log\left| \frac{1-z\bar w}{z-w}\right|
q(w)F_i(w)\, d\A(w)\,,
$$
which readily yields $\|F_{i+1}\|\le c_0 \|q\|\, \|F_i\|$.
Choosing $\varepsilon_0<\frac1{2c_0}$, we get $\|F_i\|\le
(c_0\|q\|)^i \le 2^{-i}$. Hence the series
$$
\psi = \sum_{i=0}^\infty F_i
$$
converges uniformly. Therefore $\psi$ is a weak and thus a
classical solution of the equation \eqref{eq1}. Also,
$$
\|\psi-1\|\le \sum_{i\ge 1} \|F_i\| \le
\frac{c_0\|q\|}{1-c_0\|q\|} \le 2c_0\|q\|\,.
$$
Finally,
$$
\varphi = \frac{\psi}{\|\psi\|}
$$
is the desired positive solution. \hfill $\Box$

\begin{lemma}\label{lemma2}
Let $F$ be any non-zero solution to equation \eqref{eq1}. Then
there exist a $K$-quasiconformal homeomorphism $h\colon \D\to \D$
with $h(0)=0$ and a harmonic function $U\colon \D\to \R$ such that
$F = \varphi\cdot (U\circ h)$. Moreover, the dilation $K$
satisfies
\begin{equation}\label{eq2}
K \le 1 + c_2 \|q\|\,.
\end{equation}
\end{lemma}
\par\noindent{\em Proof of Lemma~\ref{lemma2}: }
Write $F=\varphi u$, and note that (by direct computation)
equation \eqref{eq1} yields
\begin{equation}\label{eq3}
\frac{\partial}{\partial x} \left( \varphi^2 \frac{\partial
u}{\partial x}\right) + \frac{\partial}{\partial y} \left(
\varphi^2 \frac{\partial u}{\partial y}\right) = 0\,.
\end{equation}
Thus, there exists a unique smooth function $v$ with $v(0)=0$ such
that $\varphi^2 u_x = v_y$, and $\varphi^2 u_y = -v_x$. To rewrite
these equations in the complex form, we consider the
complex-valued function $w=u+iv$. An inspection shows that
$$
\frac{\partial w}{\partial\bar z} =
\frac{1-\varphi^2}{1+\varphi^2}\, \overline{\frac{\partial
w}{\partial z}}\,.
$$
In other words, $w$ satisfies the Beltrami equation
$$
\frac{\partial w}{\partial\bar z} = \mu \frac{\partial w}{\partial
z}
$$
with the Beltrami coefficient
$$
\mu = \frac{1-\varphi^2}{1+\varphi^2} \cdot
\frac{u_x+iu_y}{u_x-iu_y}\,.
$$
Clearly,
$$
|\mu| = \frac{1-\varphi^2}{1+\varphi^2} < 1\,.
$$
Since $u$ is a non-trivial solution of an elliptic equation
\eqref{eq3}, its critical points are isolated. Thus, $\mu$ is a
measurable function defined almost everywhere. By the fundamental
existence theorem \cite[Chapter~V]{Ahlfors}, there exists a
$K$-quasiconformal homeomorphism $h\colon \D \to \D$ with $h(0)=0$
such that $w=W\circ h$ where $W$ is an analytic function on $\D$.
This yields $F=\varphi\cdot(U\circ h)$ where
$U=\operatorname{Re}W$. The dilation $K$ of $h$ satisfies
$$
\frac{K-1}{K+1} \le \|\mu\|_{L^\infty}\,.
$$
Taking into account Lemma~\ref{lemma1}, we get inequality
\eqref{eq2}. \hfill $\Box$

\medskip The dilation $K$ controls geometric properties of the homeomorphism
$h$. We shall use Mori's theorem, which states that $h$ is
$\frac1{K}$-H\"older and
\begin{equation}\label{eq4}
\frac1{16} |z_1-z_2|^K \le |h(z_1)-h(z_2)| \le 16 |z_1-z_2|^{1/K}
\end{equation}
(see \cite[Section~IIIC]{Ahlfors}), and Astala's distortion
theorem \cite{Astala}:
\begin{equation}\label{eq5}
\A(h(E)) \le c_3 \A(E)^{1/K}\,.
\end{equation}
The constant $c_3$ in Astala's theorem depends on $K$ but stays
bounded when $K$ remains bounded, so we may treat it as absolute.
\begin{lemma}\label{lemma3} We have
\begin{equation}
\label{eq3.10} \A \left( \{ F>0 \} \right) \ge \frac{c_4}{\left(
\log  \beta^*(F) \right)^{1+c_2\|q\|}}\,.
\end{equation}
\end{lemma}

Later we will show how this estimate can be improved by simple
rescaling.

\medskip\par\noindent{\em Proof of Lemma~\ref{lemma3}:} We have
$\{F>0\} = h^{-1}\{U>0\}$ where $U$ is the harmonic function
obtained in the previous lemma. Hence, by the area distortion
theorem \eqref{eq5},
$$
\A \left( \{ F>0 \} \right) = \A \left( \{ h^{-1} \{U>0\}
\}\right) \ge c_5 \A \left( \{U>0 \} \right) ^K\,.
$$
By Theorem~\ref{thmharm},
$$
\A \left( \{ U>0\} \right) \ge \frac{c_6}{\log \beta^* (U, \D)}\,.
$$

Now, using Mori's theorem, we choose a positive integer $\ell_0$
so large that $h^{-1}(\frac12 \D) \supset 2^{-\ell_0} \D$. Then
$$
\frac{\max_\D |U|}{\max_{\frac12 \D} |U|} = \frac{\max_\D |U \circ
h|}{\max_{h^{-1}(\frac12 \D) } |U \circ h|} \le \frac{\max_\D
|U\circ h|}{\max_{2^{-\ell_0}\D} |U\circ h|}\,.
$$
The right hand side is bounded by
$$
c_7 \frac{\max_\D |F|}{\max_{2^{-\ell_0}\D} |F|} \le c_7 e^{
\ell_0 \beta (F)}\,.
$$
Hence $\beta^* (U, \D) \le c_8 \beta^* (F)$, and
$$
\A \left( \{ F>0\} \right) \ge \frac{c_9}{\left( \log \beta^* (F)
\right)^{K} }\,.
$$
Recalling estimate \eqref{eq2}, we get the desired result. \hfill
$\Box$

\medskip\noindent{\em The end of the proof of Theorem~\ref{thmschro}: }
The nodal set $L=\{F=0\}$ of the function $F$ contains the origin
and does not have closed loops (since it is homeomorphic to the
nodal set of the harmonic function $U$). Therefore, for any $r\in
(0,1]$, there are at least $c_{10}r^{-1}$ disjoint discs
$D_j\subset \D$ of radius $r$ centered at $z_j\in L$. For each
disc $D_j$, consider the function
$$
F_j(z)=F(z_j + rz)\,, \qquad z\in\D\,.
$$
It satisfies the equation
\begin{equation}\label{eq6}
\Delta F_j + q_j F_j = 0
\end{equation}
with $q_j(z)=r^2q(z_j+rz)$, $\|q_j\| \le r^2 \|q\| \le r^2
\varepsilon_0$. Applying Lemma~\ref{lemma3} to $F_j$ instead of
$F$ and taking into account that $\beta^*(F_j)\le \beta^*(F)$, we
get
$$
\A(\{F_j>0\}) \ge \frac{c_{4}}{\left(
\log\beta^*(F)\right)^{1+c_2\varepsilon_0r^2}}\,.
$$
To simplify the notation, denote $b=\log\beta^*(F)$ and
$s=c_2\varepsilon_0 r^2$, so that $\A(\{F_j>0\})\ge c_4 b^{-1-s}$.
Then
\begin{eqnarray*}
\A(\{ F>0 \}) &\ge& \sum_{j} \A(\{ F>0 \}\cap
D_j) \\
&=& r^2\, \sum_{j} \A(\{ F_j>0\} ) \\
&\ge& r^2 \cdot \frac{c_{10}}{r} \cdot  \frac{c_{4}}{b^{1+s}} =
\frac{c_{11}\sqrt{s} }{b^{1+s}}\,.
\end{eqnarray*}
The choice of the scaling parameter $r$ (and hence of $s$) is in
our hands. One readily checks that, for $\beta\geq\log 3$, the
function $s\mapsto \sqrt{s}b^{-s}$, $s\in(0,c_2\varepsilon_0]$,
attains its maximum at $s=(2\log b)^{-1}$ for large $b$ and at
$s=c_2\varepsilon_0$ for small $b$. In both cases the maximal
value of this function is $\ge c_{12} (\log b)^{-1/2}$. This
completes the proof. \hfill $\Box$


\section{The Donnelly-Fefferman estimate}
\label{sect4}

In this section, we prove the Donnelly-Fefferman estimate for the
doubling exponent:
\begin{thm}\label{thmweak}
For any metric disc $D\subset S$,
$$
\beta(D, f_\la ) \le a_1 \sqrt\la \,.
$$
\end{thm}

Our proof is based on a version of the  `Three Circles Theorem'
for solutions of the Schr\"odinger equation. Similar results are
known under various assumptions: see Landis \cite{Landis}, Agmon
\cite{Agmon}, Gerasimov \cite{Gerasimov}, Brummelhuis
\cite{Brummelhuis}, Kukavica \cite{Kukavica}.

Let $F$ be a solution to the equation
\begin{equation}\label{eq.3x}
\Delta F + q F = 0
\end{equation}
where $q$ is a smooth function in the unit disc $\D$. We no longer
assume that $q$ is small. Instead, the size of $q$ will be
controlled by the quantity
$$
N = \max_{\D} \left(|q| + \rho |q_\rho|\right)
$$
where $\rho$ is the polar radius. Denote $M(r)=\max_{r\D} |F|$.

\begin{thm}\label{thm4circles}
Let $F$ be a solution to the equation \eqref{eq.3x}. Then
\begin{equation}\label{eq.1x}
\frac{M(2s)}{M(s)} \le c_1 e^{c_2\sqrt N} \frac{M(8r)}{M(r)}\,,
\end{equation}
provided that $ 0<s\le r\le \frac18 $.
\end{thm}

\subsection{Proof of Theorem \ref{thm4circles}}
Our first aim will be to replace PDE \eqref{eq.3x} by a second
order ODE $\ddot{h} = L(t)h$ where $L(t)$ is a non-negative
unbounded operator on a Hilbert space such that $\dot{L}(t)$ is
also non-negative.

First, adding an extra variable $z$, we make the potential
non-positive. Put $v(x,y,z)= F(x,y) \cdot \cosh \ga z$ where $ \ga
=\sqrt N$. Then $\Delta v = (\ga^2 - q)v$, or, which is the same,
\begin{equation}\label{eq.6x}
v_{rr} + \frac2{r} v_r = - \frac1{r^2} \widetilde{\Delta}v +
(\ga^2-q) v
\end{equation}
where $r$ is the polar radius and $\widetilde{\Delta}$ is the
spherical part of the Laplacian.

Next, we make the logarithmic change of variable and put
\[
h(t, \theta) := e^{t/2} v(e^t x_\theta, e^t y_\theta, e^t
z_\theta)
\]
where $(x_\theta, y_\theta, z_\theta)=\theta\in \mathbb S^2$ and
$t\in (-\infty, 0]$. Define $Q(t, \theta) = q(e^tx_\theta,
e^ty_\theta)$. Then equation \eqref{eq.6x} turns into
\begin{equation}\label{eq.7x}
\ddot{h} = \left( - \widetilde{\Delta} + e^{2t} (\ga^2 - Q) +
\frac14 \right) h =: L(t)h\,.
\end{equation}
Note that $L(t)$ is a symmetric positive (unbounded) operator on
the Hilbert space $L^2(\mathbb S^2)$. The initial conditions for
ODE \eqref{eq.7x} are
\begin{eqnarray*}
h(-\infty) &=& 0\,, \\
\dot{h}(-\infty) &=& \lim_{r\to 0} r^{3/2} v_r = 0\,.
\end{eqnarray*}
Note that, due to our choice of $\gamma$, the derivative
$\dot{L}(t) = 2e^{2t} \left( \ga^2 - Q - \frac12 \dot{Q}\right)$
of $L(t)$ is also a non-negative operator.

At this point we make a break in the proof of the theorem and
prove a lemma on second order ODEs (cf. Agmon \cite{Agmon}):

\begin{lemma}\label{lem.ODE} Let $h$ be a solution to the equation
\[
\ddot{h} = L(t)h\,, \qquad -\infty<t\le 0\,,
\]
with
\[
h(-\infty) = \dot{h}(-\infty) = 0
\]
where $L(t)$ is a non-negative linear operator on a Hilbert space
$\mathcal H$ such that $\dot{L}(t)$ is also non-negative. Then the
function
$$
t\mapsto \log\frac{\|h\|^2}2
$$
is convex.
\end{lemma}

\noindent{\em Proof:} Denote $a(t)= \frac12 \|h\|^2$. Then
$\dot{a}(t) = (h, \dot{h})$, and
\[
\ddot{a}(t) = (h, \ddot{h}) + \|\dot h\|^2 = (h, L(t)h) +
\|\dot{h}\|^2 \ge 0\,.
\]
We need to show that $(\log a)^{\bf{\cdot\cdot}} \ge 0$, or,
equivalently, that $\ddot{a} a - \dot{a}^2 \ge 0$. We have
\begin{eqnarray*}
\ddot{a}a - \dot{a}^2 &=& \left( (L(t)h,h) + \|\dot{h}\|^2\right)
\frac{\|h\|^2}2 - (h, \dot{h})^2 \\
&\ge& \left( (L(t)h,h) + \|\dot{h}\|^2\right)
\frac{\|h\|^2}2 - \|h\|^2 \cdot \|\dot{h}\|^2 \\
&=& \left( (L(t)h,h) - \|\dot{h}\|^2\right) \frac{\|h\|^2}2\,.
\end{eqnarray*}
Further, since
\[
\frac{d}{dt} \left( (L(t)h,h) - \|\dot{h}\|^2 \right) =
(\dot{L}(t)h,h) + (L(t)\dot{h}, h) + (L(t)h, \dot{h}) -
2(\ddot{h}, \dot{h}) = (\dot{L}(t)h, h),
\]
we obtain
\[
(L(t)h,h) - \|\dot{h}\|^2 \geq \int_{-\infty}^t (\dot{L}(\tau)h,
h)\, d\tau \ge 0\,,
\]
and, thereby, $\ddot{a} a - \dot{a}^2 \ge 0$, proving the lemma.
\hfill $\Box$

\medskip\par\noindent{\em Continuation of the proof of
Theorem~\ref{thm4circles}: } Consider the spherical integral
\[
I(t) = \frac12 \iint_{\mathbb S^2} h^2(t, \theta) \, d\sigma
(\theta)
\]
($d\sigma$ is the spherical area form) for the function $h$
defined above. Since our function $h$ is even in $z$-variable, we
integrate only over the upper hemisphere $\mathbb S^2_+$.
Introduce the coordinates
\[
x=\rho \cos\varphi, \quad y=\rho \sin\varphi, \quad
z=\sqrt{1-\rho^2}
\]
on $\mathbb S^2_+$. Then
\[
d\sigma = \frac{\rho}{\sqrt{1-\rho^2}}\,d\rho d\varphi\,,
\]
and
\[
h^2(t, \theta) = e^t F^2(e^t\rho, \varphi) \cosh^2(\ga e^t
\sqrt{1-\rho^2})\,.
\]
We obtain
\begin{eqnarray*}
I(t) &=& \int_0^{2\pi} d\varphi \int_0^1
\frac{\rho\,d\rho}{\sqrt{1-\rho^2}} e^t F^2(e^t\rho, \varphi)
\cosh^2(\ga e^t \sqrt{1-\rho^2}) \\
&=& \int_0^{e^t} \frac{s\,ds}{\sqrt{e^{2t}-s^2}} \cosh^2(\ga
\sqrt{e^{2t}-s^2}) \int_0^{2\pi} F^2(s, \varphi)\, d\varphi\,.
\end{eqnarray*}
Finally, we introduce the function
\[
J(r) = I(\log r) = \int_0^r \frac{\cosh^2(\ga \sqrt{r^2-s^2})s}
{\sqrt{r^2-s^2}} \left( \int_0^{2\pi} F^2(s, \varphi)\, d\varphi
\right) ds \,,
\]
By Lemma~\ref{lem.ODE}, the function $t\mapsto \log J(e^t)$ is
convex. Hence
\begin{equation}\label{eq.9x}
\frac{J(2s)}{J(s)} \le \frac{J(2r)}{J(r)}
\end{equation}
for $0<s<r<\frac12$.

It remains to rewrite this estimate in terms of $M(r)$. For this,
we use the following standard lemma from the elliptic theory.
\begin{lemma}\label{claim}
\begin{equation}\label{eq.10x}
c_3 e^{-\sqrt{N}  r} \sqrt{\frac{J(r)}{r}} \le M(r) \le c_4 N
\sqrt{\frac{J(2r)}{2r}}\,, \qquad 0<r\le \frac12\,.
\end{equation}
\end{lemma}

\noindent{\em Proof of the upper bound:} Observe that
\begin{equation}\label{eq11x}
M(r) \le \frac{c_5 \max_{2r \D} |q|}{r} \left( \iint_{2r\D} F^2\,
d\A \right)^{1/2}\,.
\end{equation}
Indeed, after rescaling, \eqref{eq11x} reduces to its special case
when $r=\frac12$:
\begin{equation}\label{eq11y}
M(\tfrac12) \le c_6 \max_\D |q|\, \left( \iint_{\D} F^2\, d\A
\right)^{1/2}\,.
\end{equation}
To get this estimate, we represent the function $F$ as the sum of
Green's potential and the Poisson integral:
\[
F(z) =  \iint_{\rho \D} q(\zeta) F(\zeta) \log\left|
\frac{\rho^2-z\overline{\zeta}}{\rho(z-\zeta)}\right|\, d\A(\zeta)
+ \int_{\rho \mathbb T} F(\zeta)
\frac{\rho^2-|z|^2}{|\zeta-z|^2}\, dm(\zeta)\,.
\]
Here $|z|\le \frac12$, $\frac23 \le \rho \le 1$, and $m$ is the
normalized Lebesgue measure on the circle $\rho \T$. Then
$$
M(\tfrac12) \le c_7 \left(
 \max_\D |q| \iint_{\rho \D} F^2\, d\A + \int_{\rho\T} F^2\, dm
\right)\,.
$$
Averaging this by $\rho$ over $[\frac23; 1]$, we get
\eqref{eq11y}.

By definition of the function $J$,
$$
J(2r) \ge \frac1{2r}\, \iint_{2r\D} F^2\, d\A\,,
$$
and the upper bound in \eqref{eq.10x} follows from \eqref{eq11x}.

\smallskip\noindent{\em Proof of the lower bound:}
\begin{eqnarray*}
J(r) &\le& \int_0^r \frac{\cosh^2(\ga
\sqrt{r^2-s^2})}{\sqrt{r^2-s^2}}
s\,ds \cdot 2\pi M^2(r) \\
&\le& r\cosh^2 \ga r \cdot \underbrace{\int_0^1 \frac{s\,
ds}{\sqrt{1-s^2}} }_{=1} \cdot 2\pi M^2(r) \\
&\le& r e^{2\ga r} \cdot 2\pi M^2(r)\,.
\end{eqnarray*}
It remains to recall that $\ga =\sqrt{N}$. The lemma is proved.
\hfill $\Box$

\medskip\par\noindent{\em End of the proof of Theorem~\ref{thm4circles}: }
By Lemma~\ref{claim},
\[
\frac{M(2s)}{M(s)} \le c_8 N e^{\sqrt{N} s}\cdot \left(
\frac{J(4s)}{J(s)} \right)^{1/2},
\]
\[
\left( \frac{J(8r)}{J(2r)} \right)^{1/2} \le c_9 N e^{8\sqrt{N} r}
\cdot \frac{M(8r)}{M(r)}\,,
\]
and, by \eqref{eq.9x},
\[
\left( \frac{J(4s)}{J(s)} \right)^{1/2} \le \left(
\frac{J(8r)}{J(2r)} \right)^{1/2}\,.
\]
Juxtaposing these three inequalities, we obtain \eqref{eq.1x}.
\hfill $\Box$

\subsection{Proof of Theorem \ref{thmweak}}

Here and in the next section, we use the classical local
description of smooth Riemannian metrics on closed surfaces. Given
a point $p \in S$, one can choose local coordinates $(x,y), \; x^2
+y^2 \leq 1000 $, near $p$ so that the metric $g$ in these
coordinates is conformally Euclidean: $g = q(x,y)(dx^2 + dy^2)$.
The point $p$ corresponds to the origin: $p = (0,0)$. This choice
can be made in such a way that the function $q(x,y)$ is pinched
between two positive constants that depend only on metric $g$:
$$
0< q_- \leq q(x,y) \leq q_+\,,
$$
and that the $C^1$-norm of $q$ is bounded by a constant depending
only on the metric $g$.

Till the end of this section, by disc $D(p, r)$ centered at a
point $p\in S$ with radius $r$ we mean the set $\{x^2+y^2\le
r^2\}$, where $(x,y)$ are local conformal coordinates near $p$. We
can also choose our conformal charts in such way that for some
$\eta>0$ and any $p, p'\in S$ such that $\operatorname{dist}(p,
p') < \eta$, we have $D(p, \frac12) \subset D(p', 1)$.

We will refer to $(x,y)$ as {\it preferred local conformal
coordinates near $p$}. Note that in local conformal coordinates
the eigenfunction $f := f_{\lambda}$ satisfies the equation
\begin{equation}\label{eq3.2}
\Delta f + \lambda q(x,y) f = 0\,.
\end{equation}

The proof of Theorem~\ref{thmweak} is based on
Theorem~\ref{thm4circles} and the following lemma, which, in turn,
is also an easy consequence of Theorem~\ref{thm4circles} and
compactness of the surface $S$.
\begin{lemma}\label{lemmaA} For every point $p\in S$,
$$
\max_{D(p, 1)} |f_\la| \ge e^{-a_2\sqrt\la} \max_S |f_\la|\,.
$$
\end{lemma}

\par\noindent{\em Proof of Lemma~\ref{lemmaA}: } Normalize the
eigenfunction $f_\la$ by the condition
$$
\max_S |f_\la| = 1\,.
$$
Let $p_0$ be the maximum point of $|f_\la|$ on $S$. For arbitrary
$p\in S$, consider the chain of $k+1$ discs $D(p_j, 1)$ connecting
$p_0$ with $p=p_k$ in such a way that
$$
D(p_j, \tfrac1{2} ) \subset D(p_{j+1}, 1)\,, \qquad 0\le j\le
k-1\,.
$$
The number $k$ depends only on the metric $g$. Due to
Theorem~\ref{thm4circles} applied to solutions of equation
\eqref{eq3.2} with $N \le a_3 \la$,
$$
\frac{\max_{D(p_j, 1)} |f_\la|}{\max_{D(p_j, \frac12 )} |f_\la|}
\le e^{a_4 \sqrt \la} \frac{\max_{D(p_j, 8)} |f_\la|}{\max_{D(p_j,
1)} |f_\la|} \le \frac{e^{a_4 \sqrt\la} }{\max_{D(p_j, 1)}
|f_\la|}\,,
$$
or
$$
\left( \max_{D(p_j, 1)} |f_\la| \right)^2 e^{-a_4\sqrt \la} \le
\max_{D(p_j, \frac12 )} |f_\la| \le \max_{D(p_{j+1}, 1)}
|f_\la|\,.
$$
Making $k$ iterations, we arrive at
$$
e^{-a_2\sqrt \la} \le \max_{D(p, 1)} |f_\la|\,,
$$
proving the lemma. \hfill $\Box$

\medskip\par\noindent{\em Proof of Theorem~\ref{thmweak}: }
Fix a disc $D(p, s)\subset S$ of radius $s$ centered at $p$. Since
in each conformal chart the Riemannian metric is equivalent to the
Euclidean one: $q_- (dx^2 + dy^2) \le g \le q_+ (dx^2 + dy^2)$, it
suffices to show that
\begin{equation}\label{eq3b}
\frac{\max_{D(p, s)} |f_\la|}{\max_{D(p, \frac12 s)} |f_\la| } \le
e^{a_5\sqrt\la}
\end{equation}
provided that $\la \ge 2$. The previous Lemma yields \eqref{eq3b}
for $s \ge 2 $. Assume now that $s<2$. We apply
Theorem~\ref{thm4circles} to the solution $f_\lambda$ of equation
\eqref{eq3.2} in the disk $D(p,16)$. As above, the parameter $N$
in Theorem~\ref{thm4circles} does not exceed $a_3\la$, and we
immediately get \eqref{eq3b}. \hfill $\Box$


\section{Proof of Theorem~\ref{thmsymmet}}\label{sect5}

It suffices to prove Theorem~\ref{thmsymmet} assuming that $D$ is
a Euclidean disc lying in a chart with preferred local coordinates
and that the center of $D$ belongs to the nodal line. All metric
notions (distance, area and discs) pertain to the Euclidean
metric.

We fix a positive number $\rho_0$ depending only on metric $g$
such that
\begin{equation}\label{eq4*}
\rho_0^2 q_+ \le \varepsilon_0
\end{equation}
where $\varepsilon_0$ is the numerical constant that controls
`smallness' of the potential in Section~\ref{sect3}.
\begin{definition}\label{def} {\rm A {\em good} disc on
$S$ is a disc of radius $\le \rho_0 \lambda^{-1/2}$ whose center
lies on the nodal line $\{f_\lambda=0\}$.}
\end{definition}

\begin{lemma}\label{lemmagood}
Let $D$ be a good disc. Then
\begin{equation}\label{eq3.1}
\frac{\A (S_+(\la)\cap D)}{\A (D)} \ge \frac{a_1}{\log\la \cdot
\sqrt{\log\log\la } }
\end{equation}
\end{lemma}

\par\noindent{\em Proof of Lemma~\ref{lemmagood}: }
Given a good disc $D$ of radius $r\la^{-1/2}$, $r\in (0, \rho_0)$,
with center $p$, define the function $F$ on the unit disc
$\D\subset\C$ by
\[
F(x,y) = f_\la (\tfrac{rx}{\sqrt{\la}}, \tfrac{ry}{\sqrt{\la}})\,,
\]
where $(x,y)$ are preferred local coordinates near $p$. Then
$$
\Delta F + r^2 q ( \tfrac{rx}{\sqrt{\la}}, \tfrac{ry}{\sqrt{\la}})
F = 0\,,
$$
and $F(0,0)=0$. We have
$$
\frac{\A(S_+(\la)\cap D)}{\A(D)} \ge a_2\,
\frac{\A(\{F>0\})}{\A(\D)}\,.
$$
Due to the choice of $\rho_0$, we can apply
Theorem~\ref{thmschro}:
$$
\frac{\A(\{F>0\})}{\A(\D)} \ge  \frac{a_3}{\log\beta^*(F) \cdot
\sqrt{ \log\log\beta^*(F)} }\,.
$$
It follows from Theorem~\ref{thmweak} that the right hand side is
$$
\ge \frac{a_4}{\log\la \cdot \sqrt{\log\log\la} }\,.
$$
This proves the lemma. \hfill $\Box$

\medskip Now we are ready to prove Theorem~\ref{thmsymmet}.

\medskip\par\noindent{\em Proof of Theorem~\ref{thmsymmet}: }
In the proof we use the fact that the inradius of every nodal
domain does not exceed $\rho_1 \la^{-1/2}$ where $\rho_1$ depends
only on metric $g$ (see \cite{Bruning}). Let $D$ be a disc of
radius $R$ centered at the nodal line.

\smallskip\par\noindent{\em Case I: $R\ge 100\rho_1 \la^{-1/2}$. }
Consider the collection of all discs of radii $2\rho_1 \la^{-1/2}$
with centers on the nodal line $L=\{f_\la=0\}\cap\frac12 D$. Let
$\left\{D_i\right\}_{i=1,...,N}$ be a maximal subcollection of
pairwise disjoint discs. We claim that every point $p\in\frac14D$
lies at distance at most $6\rho_1 \la^{-1/2}$ from the center of
some $D_i$. Indeed, otherwise, choose a point $p'\in L$ with
$\hbox{dist}(p,p')\le \rho_1\la^{-1/2}$, and consider the disc
$D'$ of radius $2\rho_1\la^{-1/2}$ centered at $p'$. Our
assumption yields that $D'$ is disjoint from all $D_i$'s, which
contradicts the maximality of the subcollection. The claim
follows.

The claim yields that the discs $\left\{ 4D_i\right\}$ cover
$\frac14 D$, so we get the inequality
\begin{equation}\label{eq3.3}
\sum_i \A(4D_i) \ge \A( \tfrac14 D)\,.
\end{equation}
Denote by $D'_i$ the good disc $D'_i = \frac{\rho_0}{2\rho_1}\,
D_i$. Note that $\A (D'_i)\ge a_5 \A(4D_i)$, and $\A(\frac14 D)
\ge a_6 \A(D)$. Therefore, using \eqref{eq3.3}, we get
\begin{equation}\label{eq3.4}
\sum_i \A(D'_i) \ge a_7 \A(D)\,.
\end{equation}
Further,
$$
\A(S_+(\la)\cap D) \ge \sum_i \A(S_+(\la) \cap D'_i) \ge
\frac{a_8}{\log\la \cdot \sqrt{\log\log\la} } \sum_i \A(D_i')\,,
$$
where the last inequality follows from Lemma~\ref{lemmagood}.
Combining this with \eqref{eq3.4}, we obtain
$$
\A(S_+(\la)\cap D) \ge \frac{a_9}{\log\la \cdot \sqrt{
\log\log\la} }\, \A(D)\,,
$$
which proves the theorem in this case.

\smallskip\par\noindent{\em Case II: $R\le 100\rho_1 \la^{-1/2}$. }
Choose $\rho_0$ in Definition~\ref{def} of good discs to be less
than $400\rho_1$. Then the disc $D'$ concentric with $D$ of radius
$$
r = R \cdot \frac{\rho_0}{400\rho_1}
$$
is good. Applying Lemma~\ref{lemmagood}, we get
\begin{eqnarray*}
\A\left( S_+(\la)\cap D\right) &\ge& \A\left( S_+(\la) \cap D'\right) \\
&\stackrel{\eqref{eq3.1}}\ge& \frac{a_1 \, \A(D')}{\log\la \cdot
\sqrt{\log\log\la} } \ge \frac{a_{10}\, \A(D)}{\log\la \cdot
\sqrt{\log\log\la} }\,.
\end{eqnarray*}
as required. This completes the proof in Case~II, finishing off
the proof of Theorem~\ref{thmsymmet}. \hfill $\Box$


\section{Logarithmic asymmetry}\label{sect6}

In this section, we prove the results confirming sharpeness of our
lower bounds for the area of positivity. First, we shall construct
harmonic polynomials with small positivity area:
\begin{thm}\label{thmI}
There exists a sequence of complex polynomials $P_N(z)$, $N=2$,
$3$, ..., such that $\hbox{deg}P_N=N$, $P_N(0)=0$, and
$$
\A(\{\operatorname{Re}P_N>0\}\cap \D) \le \frac{c}{\log N}\,.
$$
\end{thm}

Of course, Theorem~\ref{thmI} yields the upper bound for the
Nadirashvili  constant $\mathcal N$ in Theorem~\ref{thmnadconst}.
Then we prove Theorem~\ref{thmsharp} `transplanting' the
polynomials $P_N$ to a small chart on the sphere $\mathbb S^2$ and
transforming  them into spherical harmonics on $\mathbb S^2$.

\subsection{Proof of Theorem~\ref{thmI} }
Let us explain the idea behind the construction of harmonic
polynomials in Theorem~\ref{thmI}. We start with an entire
function $E(z)$ on $\C$ which is bounded outside a semi-strip
$\Pi_+=\{x\ge 0, \ |y|\le \frac{\pi}2\}$. For simplicity, assume
that, for all sufficiently large $R$, the maximum $M(R) =
\max_{R\D} |E|$ is attained at $z=R$ and $E(R)$ is real positive.
Fix a sufficiently large $R$, and note that the function $G(z) :=
E(z+R) - E(R)$ vanishes at $0$ and $\mbox{Re} G<0$ outside the
strip $\Pi = \{ |y|\le \frac{\pi}2 \}$. We will check that $G$
admits a good approximation on the disc $R\D$ by its Taylor
polynomial $Q_N$ of degree $N \approx \log M(R)$, so the set
$\{\mbox{Re} Q_N>0 \}\cap R\D$ is still contained in the strip
$\Pi$. Rescale the polynomial $Q_N$ and set $P_N(z) := Q_N(Rz)$.
Then
\begin{eqnarray}\label{eq6.1a}
\frac{\A(\{\mbox{Re} P_N>0\}\cap \D)}{\A(\D)} &=&
\frac{\A(\{\mbox{Re} Q_N>0\}\cap R\D)}{\A(R\D)} \nonumber \\
&\le& \frac{\A(\Pi \cap R\D)}{\A(R\D)} \le \frac{c_1}{R} \approx
\frac{c_2}{M^{-1}(e^N)}
\end{eqnarray}
where $M^{-1}$ is the inverse function to the function $M$.

To get the optimal example, we have to minimize the right hand
side of \eqref{eq6.1a}, that is, to start with the function $E$ as
above with the minimal possible growth. According to the
Phragm\'en-Lindel\"of principle, the minimal growth rate for
$M(R)$ is of the double exponent order $\exp\exp R$. Therefore,
$N\approx \exp R$, and \eqref{eq6.1a} yields
$$
\frac{\A(\{\mbox{Re} P_N>0\}\cap \D)}{\A(\D)} \le \frac{c_3}{\log
N}\,,
$$
as needed.

Now, let us pass to the formal construction.

\medskip\par\noindent{\em Proof of Theorem~\ref{thmI}: }
Following Mittag-Leffler and Malmquist, we produce an entire
function $E(z)$ such that
\begin{equation}\label{eq6.2}
|E(z)| \le c_4 \qquad \mbox{for} \quad z\notin \Pi_+\,,
\end{equation}
and
\begin{equation}\label{eq6.3}
\left| E(z) - e^{e^z} \right| \le c_4 \qquad \mbox{for} \quad
z\in\Pi_+\,.
\end{equation}
To get $E$, denote $\Pi' = \{x>0, \ |y|\le \frac23 \pi \}$, $\Pi''
= \{x>-1, \ |y|\le \frac43 \pi \}$, and consider a smooth cut-off
function $\chi$ on $\C$ that equals $1$ on $\Pi'$ and vanishes
outside $\Pi''$. In addition, choose $\chi$ so that $\left|
\bar\partial \chi \right|$ is uniformly bounded. Define
\begin{equation}\label{eq6.4}
u(z) = \frac1{\pi} \iint_{\C} \frac{e^{e^\zeta} \bar\partial \chi
(\zeta)}{z-\zeta}\, d\A(\zeta) = \frac1{\pi} \iint_{\Pi''
\setminus \Pi'} \frac{e^{e^\zeta} \bar\partial \chi
(\zeta)}{z-\zeta}\, d\A(\zeta)\,.
\end{equation}
Then
$$
\bar\partial u = e^{e^z} \bar\partial\chi = \bar\partial\left(\chi
e^{e^z}\right)\,,
$$
so the function $E(z)=\chi \exp\exp z - u(z)$ is entire. To
establish properties \eqref{eq6.2} and \eqref{eq6.3}, it suffices
to show that $|u|$ is bounded. This readily follows from the fact
that $|\exp\exp z| = \exp(e^x \cos y)$, and, therefore, the
integrand in \eqref{eq6.4} decays very rapidly when
$\zeta\to\infty$ within the layer $\Pi''\setminus \Pi'$.

Now choose $R>0$ such that $R> 2c_4 + 1$, and set $G(z) =
E(z+R)-E(R)$. Then $G(0)=0$, and, for $z\notin\Pi :=\{|
\operatorname{Im} z|\le \frac{\pi}2 \}$,
$$
\mbox{Re} G(z) \le c_4 - (e^{e^R}-c_4) \le -1\,.
$$
If $r$ is sufficiently large, we have
$$
M(r) := \max_{r\D} |G(z)| \le \exp (c_5 e^r)\,.
$$

It remains to approximate $G$ by its Taylor
polynomial\footnote{This step is not needed for the upper bound
for the Nadirashvili  constant $\mathcal N$ in
Theorem~\ref{thmnadconst} that can be obtained directly by scaling
$\mbox{Re} G$.}. Let
$$
G(z) = \sum_{n=1}^\infty a_n z^n\,, \qquad Q_N(z) = \sum_{n=1}^N
a_n z^n\,, \qquad R_N(z) = G(z) - Q_N(z)\,.
$$
By Cauchy's inequalities,
$$
|a_n| \le \frac{M(\rho)}{\rho^n} \le \exp (c_5 e^\rho - n\log
\rho)\,.
$$
Assume that $n$ is sufficiently large, and choose $\rho$ so that
$c_5 \rho e^\rho=n$. Then
$$
\frac{n}{c_5\log n} \le e^\rho \le n\,,
$$
and
$$
|a_n| \le \exp \left( c_5 n - n \log\log \frac{n}{c_5 \log
n}\right) \le  \left( \frac{c_6}{\log n}\right)^n \,.
$$
Therefore, for $|z|=r$,
$$
|R_N(z)| \le \sum_{n=N+1}^\infty |a_n| r^n \le \sum_{n=N+1}^\infty
\left( \frac{c_6 r}{\log n} \right)^n\,.
$$
Hence, if $r\le r_N := \frac12 c_6^{-1} \log N$ and $N$ is large
enough, we have $|R_N(z)|\le \frac12$, which yields $\mbox{Re}
Q_N(z) \le - \frac12$ for $|z|\le r_N$, $|\mbox{Im} z| \ge
\frac{\pi}2$. Finally, make a rescaling $P_N(z) = Q_N(r_N z)$.
This is a polynomial of degree $N$ with $P_N(0) = Q_N(0) = 0$, and
\begin{eqnarray*}
\frac{\A(\{\mbox{Re} P_N> - \frac12\} \cap\D)}{\A (\D)} &=&
\frac{\A(\{\mbox{Re} Q_N> - \frac12\}\cap r_N\D)}{\A (r_N\D)} \\
&\le& \frac{\A(\Pi\cap r_N\D)}{\A(r_N \D)} \le \frac{2\pi r_N}{\pi
r_N^2} = \frac{2}{r_N} = \frac{2c_6}{\log N}\,.
\end{eqnarray*}
This completes the proof of Theorem~\ref{thmI}. \hfill $\Box$

\subsection{Proof of Theorem~\ref{thmsharp} }
We work on the sphere $\mathbb S^2 = \{ x_1^2 + x_2^2 + x_3^2 =
1\}$ endowed with the standard spherical metric. The spectrum of
the Laplacian on $\mathbb S^2$ is given by $\la_N = N(N+1)$, where
each $\la_N$ has multiplicity $2N+1$. Put $A = (0, 0, 1)$, and
consider the upper hemi-sphere $\mathbb S^2_+=\{x_3>0\}$. Then
$(x_1, x_2)$ are local coordinates on $\mathbb S^2_+$. Put $z =
x_1 + ix_2$, $r = \sqrt{x_1^2 + x_2^2}$, $z = re^{i\theta}$.
Consider the space $\mathcal S_N$ of complex valued spherical
harmonics corresponding to the eigenvalue $\la_N$ that vanish at
$A$. Clearly, $\hbox{dim}_\C\mathcal S_N = 2N$. We shall use the
following classical

\begin{lemma}\label{lemma6.1}
There exists a basis $e_1$, $e_2$, ...., $e_N$, $e_{-1}$,
$e_{-2}$, ..., $e_{-N}$ in $\mathcal S_N$ such that each function
$e_j$ restricted to $\mathbb S^2_+$ has the form
$$
e_{j}(x_1, x_2) = L_N^{(j)}(\sqrt{1-r^2}) z^j\,, \qquad j=1, ...,
N\,,
$$
$$
e_{-j}(x_1, x_2) = L_N^{(j)}(\sqrt{1-r^2}) \bar{z}^j\,, \qquad
j=1, ..., N\,.
$$
Here $L_N$ is the Legendre polynomial of degree $N$, and
$L_N^{(j)}$ stands for its $j$-th derivative.
\end{lemma}
For the proof, see, e.g., \cite[Lemma~3.5.3]{Groemer}.

\medskip\par\noindent{\em Proof of
Theorem~\ref{thmsharp}: } Write $L_N^{(j)}(\sqrt{1-r^2}) = A_{jN}
+ rB_{jN}(r)$, where $B_{jN}$ is a continuous function on $[0;
1]$, and $A_{jN} = L_N^{(j)}(1)$.  Since all zeroes of the
Legendre polynomial $L_N$ are real and lie in the interval
$(-1;1)$, and since its leading coefficient is positive, we have
$A_{jN}>0$. In view of Theorem~\ref{thmI}, there exist a sequence
of complex polynomials $P_N(z) = \sum_{j=1}^N \alpha_{jN} z^j$ and
a sequence of small positive values $\{\kappa_N\}$ such that
\begin{equation}\label{eq6.1}
\A(\{z\in\D\colon \mbox{Re} P_N(z) > - \kappa_N\}) \le
\frac{c_7}{\log N}\,.
\end{equation}
Let $\delta=\delta_N$ be a sufficiently small positive number (to
be chosen later). Fix $N$ large enough, and consider the spherical
harmonic $f_N\in\mathcal S_N$ defined by,
$$
f_N(r, \theta) = \sum_{j=1}^N \beta_j L_N^{(j)}(\sqrt{1-r^2}) r^j
e^{ij\theta}
$$
with
$$
\beta_j = \frac{\alpha_{j}}{A_j \delta^j}
$$
(to simplify notation, we suppress the subindex $N$ for the
coefficients $\alpha_j$, $\beta_j$ and $A_j$, as well as for the
functions $B_j$). Rescaling, define
$$
F_N(r, \theta) = f_N(\delta r, \theta) = \sum_{j=1}^N \beta_{j}
\left( A_j + \delta r B_j( \delta r) \right) \delta^j r^j
e^{ij\theta}\,.
$$
Then
$$
\left| F_N(r, \theta) - P_N\left(re^{i\theta}\right)\right| =
\left| \sum_{j=1}^N \frac{\alpha_{j} \delta r B_j(\delta r)}{A_j}
r^j e^{ij\theta} \right| \le \delta M_N
$$
where
$$
M_N = \max_{1\le j \le N} \max_{r\in [0; 1]} \left|
\frac{B_j(r)}{A_j}\right| \cdot \sum_{j=1}^N |\alpha_{j}|\,.
$$

Choose $\delta=\delta_N$ so small that
$$
\delta M_N < \kappa_N\,.
$$
Put $E_N = \{|z|<\delta_N\}$, $D_N = \{(x_1, x_2, x_3)\in \mathbb
S^2 \colon x_1+ ix_2 \in E_N \}$. We have (writing $\A_e$ and
$\A_s$ for the euclidean and the spherical areas respectively)
\begin{eqnarray*}
\frac{\A_e (\left\{ \mbox{Re} f_N>0 \right\}\cap E_N)}{\A_e(E_N)}
&=& \frac{\A_e (\left\{ \mbox{Re} F_N>0 \right\}\cap \D)}{\A_e(\D)} \\
&\le& \frac{\A_e (\left\{ \mbox{Re} P_N > -\kappa_N \right\}\cap
\D)}{\A_e(\D)} \le  \frac{c_7}{\pi \log N}\,.
\end{eqnarray*}
At the same time,
$$
\frac{\A_e (\left\{ \mbox{Re} f_N>0 \right\}\cap E_N)}{\A_e(E_N)}
\ge c_8 \cdot \frac{\A_s (\left\{ \mbox{Re} f_N>0 \right\}\cap
D_N)}{\A_s(D_N)}\,.
$$
This yields
$$
\frac{\A_s (\left\{ \mbox{Re} f_N>0 \right\}\cap D_N)}{\A_s(D_N)}
\le \frac{c_7 \cdot c_8^{-1}}{\pi \log N}
$$
for all $N$, as required. \hfill $\Box$


\section{Discussion and questions}\label{sect7}

\subsection{Quasi-conformal or $C^1$-smooth?}
\label{subsec7.1}

The link between harmonic functions and Laplace-Beltrami
eigenfunctions on surfaces given in Lemma~\ref{lemma2} above plays
a crucial role in the present paper. Recall that the lemma states
that, for any solution $F$ to the Schr\"odinger equation $\Delta F
+ qF = 0$ in the disc $\D$ with small smooth potential $q$, there
exist a harmonic function $U\colon \D\to\R$, a positive function
$\varphi$, and a quasi-conformal homeomorphism $h$ of $\D$ such
that $F = \varphi \cdot (U \circ h)$. It remains unclear  to us
whether $h$ can be chosen to be $C^1$-smooth (or Lipschitz) with
controlled differential:
$$
\|dh\|, \|dh^{-1}\| \le C(\|q\|).
$$
Such a result would immediately remove the double logarithm in the
area estimates presented in Theorems~\ref{thmsymmet} and
\ref{thmschro}. The refined estimates would be sharp in view of
Theorem~\ref{thmsharp}.


\subsection{At which scale does quasisymmetry break?}

Recall that Theorem \ref{thmsharp} establishes the existence of a
sequence of spherical harmonics $\{f_i\}$ on the 2-sphere
corresponding to eigenvalues $\lambda_i \to \infty$ and a sequence
of discs $D_i \subset \mathbb{S}^2$ such that each $f_i$ vanishes
at the center of $D_i$ and
$$ \frac{\A(S_+(\lambda_i) \cap D_i)}{\A(D_i)} \leq
\frac{C}{\log \lambda_i}\,.
$$
In our proof, the radii $r_i$ of the discs $D_i$ decay  very
rapidly as the eigenvalues $\lambda_i$ tend to infinity. It would
be interesting to explore what is the optimal (that is, the
slowest) possible rate of decay of the $r_i$'s for which the
inequality above is still valid. For instance, can this happen on
the wave-length scale $r_i \sim 1/\sqrt{\lambda_i}\;$? Let us
emphasize that the sequence $r_i$ must converge to zero in view of
the fact that spherical harmonics enjoy quasisymmetry
$$ \frac{\A(S_+(\lambda_i) \cap D)}{\A(D)} \geq
\text{const}(r)\,
$$
for any disc $D$ of radius $\geq r$ (this follows from \cite{DF1}
and \cite{N} since the spherical metric is real analytic).

\subsection{The doubling exponent: from uniform measurements to
statistics}

The next discussion is a result of our attempt to digest
Nadirashvili's approach in \cite{N}. Start with a sequence of
eigenfunctions $f_{\lambda},\; \lambda \to +\infty$. Fix $r>0$
small enough, and consider the function
$$b(x,\lambda) := \beta(D(x, \tfrac{r}{\sqrt{\lambda}}),f_{\lambda})$$
where $D(x, \frac{r}{\sqrt{\lambda}})$ stands for the metric disc
of radius $ \frac{r}{\sqrt{\lambda}}$ with the center at the point
$x \in S$. Recall the Donnelli-Fefferman estimate
$$
B_\infty(\lambda):= \sup_{x \in S} b(x,\lambda) \leq
c\sqrt{\lambda},
$$
which played a crucial role in our approach. Interestingly enough,
replacing the $L_{\infty}$-norm of $b$ by the $L_1$- norm
$$B_1 (\lambda) := \frac{1}{\A(S)}\int_S b(x,\lambda) d\A(x)\,,$$
we get a quantity that is closely related to the length of the
nodal line $L_{\lambda}:= \{f_{\lambda} = 0\}$:
\begin{equation}
\label{eqleng} C^{-1} \cdot
\text{Length}(L_{\lambda})\lambda^{-1/2} - C \leq B_1(\lambda)
\leq C \cdot \text{Length}(L_{\lambda})\lambda^{-1/2} + C,
\end{equation}
where the constant $C>1$ depends only on metric $g$ on the surface
$S$.

Here is a sketch of the proof. Define $N(x,\lambda)$ to be the
number of intersection points between the boundary circle $T$ of
the disc $D:= D(x, \frac{r}{\sqrt{\lambda}})$, and the nodal line
$L_{\lambda}$. Let $U$ be the harmonic function associated to
$f_{\lambda}|_D$ as in Lemma~\ref{lemma2}. Then $N(x,\lambda)$
equals the number of sign changes of $U$ on $\partial \D$. Using
the topological interpretation of the doubling exponent of
harmonic functions (see formula \eqref{eqgelfond}) it is possible
to show that
$$N(x,\lambda) \simeq \beta(\D,U) \simeq
\beta(D,f_{\lambda})= b(x,\lambda).$$ Applying an elementary
integral geometry argument we get
$$\text{Length}(L_{\lambda}) \simeq \int_S N(x,\lambda)
d\A(x) \cdot \sqrt{\lambda}  \simeq B_1(\lambda) \sqrt{\lambda},$$
which readily yields inequality \eqref{eqleng}.

Inequality \eqref{eqleng} clarifies the function-theoretic meaning
of the Yau conjecture for surfaces, which states that
$\text{Length}(L_{\lambda}) \simeq \sqrt{\lambda}$: {\em the
expectation of the doubling exponent of $f_\la$ on a random metric
disc of radius $\sim 1/\sqrt{\lambda}$ is bounded by a constant
depending only on the Riemannian metric.} The Yau conjecture was
proved in \cite{DF1} in any dimension for real analytic metrics
$g$.

\subsection{What happens in higher dimensions?}

It would be interesting to extend Theorem \ref{thmsymmet} and to
explore the local asymmetry of the sign distribution  for
eigenfunctions on higher-dimensional manifolds. This problem has
the following counterpart for harmonic functions on the unit ball
$\mathbb B\subset \R^n$, $n\ge 3$. Let $u$ be a non-zero harmonic
function on $\mathbb B$ vanishing at the origin. What is the
optimal bound for $\mbox{Vol} (\{u>0\})$ in terms of its doubling
exponent $\beta(\mathbb B, u)$? Using Carleman's method
\cite{Carleman} or otherwise, one can easily show that
$$
\mbox{Vol} (\{u>0\}) \ge \frac{c}{(\beta(\mathbb B, u))^{n-1}}\,.
$$
However, we believe that this estimate is very far from being
sharp.

\bigskip\par\noindent F\"edor Nazarov\\
Department of Mathematics\\
Michigan State University\\
East Lansing, MI 48824\\
USA \\
\smallskip
\texttt{\small fedja@math.msu.edu}

\bigskip\par\noindent
Leonid Polterovich\\
School of Mathematics\\
Tel Aviv University\\
Tel Aviv 69978\\
Israel \\
\smallskip
\texttt{\small polterov@post.tau.ac.il}

\bigskip\par\noindent
Mikhail Sodin\\
School of Mathematics\\
Tel Aviv University\\
Tel Aviv 69978\\
Israel\\
\smallskip
\texttt{\small sodin@post.tau.ac.il}

\end{document}